# Abstract interpolation Thiele type fraction

V.L. Makarov[1], I.I. Demkiv[2]
[1]*Instityte of mathematics NAS of Ukraine, Kyiv*
[2] *Lviv polytechnic National University, Lviv, Ukraine*
*November 21, 2015*

**Abstract.** Abstract continued Thiele type fraction has been constructed, which is an interpolation one for nonlinear operator acting from linear topological space $X$ to algebra $Y$ with a unit. In particular cases it changes into both a classic Thiele faction and a matrix-valued Thiele-type fraction from multiple variables.

## 1. Introduction.

A number of authors have been dealing with generalization of Thiele fractions; see for instance, [1-7] and others. These generalizations can be conditionally divided into two classes. The first class involves the works on generalization of Thiele fractions in case of multivariable functions, two as a rule (see [1-4, 7]). The second class includes the works on generalization Thiele fractions in case of vector-valued matrix-valued functions from one variable (see [5, 6]). Moreover, there are some results, devoted to construction of matrix-valued interpolants from two variables [8]. However, all fractional interpolants suggested in the abovementioned works, unlike the classic Thiele fraction, have a significant drawback: when replacing the last interpolation node with an arbitrary element from the relevant set of definition, the interpolant is not converted into ordinary (vector-valued, matrix-valued) function which interpolates. Note also that the task of construction of vector-valued or matrix-valued interpolants is equivalent to traditional interpolation task, and, therefore, the need to build vector-valued or matrix-valued interpolants should be justified by each particular application.

The aim of this work is a generalization of Thiele fractions in the case of interpolation of non-linear operators acting from linear topological space $X$ to algebra $Y$ with a unit i, which does not have the specified drawback. Hence, as a special case an interpolation Thiel type fraction is obtained for functions with any number of variables without geometric restrictions on the location of interpolation nodes.

## 2. Abstract interpolation Thiele type fraction.

Let us start this article with constructive considerations concerning construction of the most general construction of interpolation Thiele type continued fraction (ICF). Let us consider the "two-storey" fraction

$$T_2(u) = F(u_0) + l_1(u-u_0)[I + l_2(u-u_1)]^{-1}, \qquad (1)$$

where $l_1, l_2$ – linear, $F$ – nonlinear operators, acting from linear topological space X to algebra Y with a unit I, elements $u, u_0, u_1 \in X$. For operator $F$ its values $F(u_{i-1,i}(\xi_i)),\ \xi_i \in [0,1],\ \ i = 1,2$ on continual knots are known

$$u_{i-1,i}(\xi_i) = u_{i-1} + g_{\xi_i}(u_i - u_{i-1}),\ \xi_i \in [0,1], i = 1,2. \qquad (2)$$

Here $g_z$ – linear, differentiated by z operator, which acts from X to X, and has the properties

$$\mathrm{g}_0 = \mathrm{E},\ \ \mathrm{g}_1 = 0,\ \ \mathrm{g}_\tau \mathrm{g}_\xi = \mathrm{g}_{\max(\tau,\xi)},\ \tau, \xi \in [0,1], \qquad (3)$$

where E: X → X identity operator. Examples of operators $g_z$ with properties (3) for the case of Hilbert space X = H and the space of piecewise continuous functions Q[0, 1] see in [12, 13]. Linear operators $l_1, l_2$ are given by the formulas

$$l_1(u-u_0) = -\int_0^1 F_1'(u_0 + g_{\tau_1}(u_1-u_0))dg_{\tau_1}(u-u_0), \quad F_1(u) = F(u),$$

$$l_2(u-u_1) = -\int_0^1 F_2'(u_1 + g_{\tau_2}(u_2-u_1))dg_{\tau_2}(u-u_1),$$

$$F_2(u) = \left[F(u) - F(u_0)\right]^{-1} l_1(u-u_0) \qquad (4)$$

and define on the set of operators twice differentiated by Gateaux, for which there are integrals (4), divided differences of the first order (see [10,11]). Let us check the fulfillment of interpolation conditions. Let us substitute in (1)-(4) a continual set $u_{1,2}(\xi_2)$ and use a property (3) of operator $\mathrm{g}_\tau$, then we will get

$$l_1(u_1) = -F(u_0) + F(u_1),$$

$$l_2(u_{1,2}(\xi_2) - u_1) = -\int_0^1 F_2'(u_1 + g_{\tau_2}(u_2-u_1))dg_{\tau_2}g_{\xi_2}(u_2-u_1) =$$

$$= -\int_{\xi_2}^1 F_2'(u_1 + g_{\tau_2}(u_2-u_1))dg_{\tau_2}(u_2-u_1) =$$

$$= -\int_{\xi_2}^1 \frac{d}{d\tau_2} F_2(u_1 + g_{\tau_2}(u_2-u_1))d\tau_2 =$$

$$= -F_2(u_1) + F_2(u_{1,2}(\xi_2)) =$$

$$= -\left[F(u_1) - F(u_0)\right]^{-1} l_1(u_1 - u_0) +$$

$$+\left[F((u_{1,2}(\xi_2)) - F(u_0)\right]^{-1} l_1(u_{1,2}(\xi_2) - u_0) =$$

$$= -I + \left[F((u_{1,2}(\xi_2)) - F(u_0)\right]^{-1} l_1(u_{1,2}(\xi_2) - u_0).$$

The latter relations with (1) lead to the conclusion that

$$T_2(u_{1,2}(\xi_2)) = F(u_{1,2}(\xi_2)),\ \forall \xi_2 \in [0,1],$$

i.e. (4) is an abstract two-storey Thiele type interpolation fraction with continual interpolation knot $u_{1,2}(\xi_2)$ and the common interpolation knot $u_0$ (the latter is obvious).

In the general case of $n$-storey fraction

$$T_n(u) = F(u_0) + l_1(u - u_0)[I + l_2(u - u_1) \times$$
$$\times [I + l_3(u - u_2) ... [I + l_n(u - u_{n-1})]^{-1} ... ]^{-1}]^{-1} =$$
$$= \underset{p=1}{\overset{n}{D}} \frac{l_p\left(u - u_{p-1}\right)}{I} \tag{5}$$

(the last formula is a symbolic representation) its components are determined the following way

$$l_k(u - u_{k-1}) = -\int_0^1 F_k'(u_{k-1} + g_{\tau_k}(u_k - u_{k-1})) dg_{\tau_k}(u - u_{k-1}),$$

$$k = 1, 2, ..., n, \quad F_1(u) = F(u), \tag{6}$$

$$F_i(u) =$$
$$\left[-I + ... \left[\left[-I + l_0(u)\right]^{-1} l_1(u - u_0)\right]^{-1} ... l_{i-2}(u - u_{i-3})\right]^{-1} l_{i-1}(u - u_{i-2})$$

$$i = 1, 2, ..., \quad u_{-1} = 0, \quad l_0(u) = -F(u_0) + F(u) + I. \tag{7}$$

Then, using mathematical induction, let us suppose that

$$T_n(u_i) = T_i(u_i) = F(u_i), \quad i = 0, 1, ..., k \tag{8}$$

and prove the validity of the relation (8) given $i = k + 1$. We get

$$T_n(u_{k+1}) = T_{k+1}(u_{k+1}),$$

$$l_{k+1}(u_{k+1}-u_k)=-\int_0^1 F_{k+1}{}'(u_k+g_{\tau_{k+1}}(u_{k+1}-u_k))dg_{\tau_{k+1}}(u_{k+1}-u_k)=$$

$$=F_{k+1}(u_{k+1})-F_{k+1}(u_k)=\mathop{D}_{p=1}^{k+1}\frac{l_{k+1-p}\left(u_{k+1}-u_{k-p}\right)}{-I}-\mathop{D}_{p=1}^{k+1}\frac{l_{k+1-p}\left(u_k-u_{k-p}\right)}{-I}=$$

$$=\mathop{D}_{p=1}^{k+1}\frac{l_{k+1-p}\left(u_{k+1}-u_{k-p}\right)}{-I}-\left[-I+\mathop{D}_{p=2}^{k+1}\frac{l_{k+1-p}\left(u_k-u_{k-p}\right)}{-I}\right]^{-1} l_k\left(u_k-u_{k-1}\right)=$$

$$=\left[-I+\mathop{D}_{p=2}^{k+1}\frac{l_{k+1-p}\left(u_{k+1}-u_{k-p}\right)}{-I}\right]^{-1} l_k\left(u_{k+1}-u_{k-1}\right)-I\,.$$

Further,

$$I+l_k(u_{k+1}-u_{k-1})\left[I+l_{k+1}(u_{k+1}-u_k)\right]^{-1}=\mathop{D}_{p=2}^{k+1}\frac{l_{k+1-p}\left(u_{k+1}-u_{k-p}\right)}{-I}=$$

$$=\mathop{D}_{p=1}^{k}\frac{l_{k-p}\left(u_{k+1}-u_{k-p-1}\right)}{-I},$$

$$I+l_{k-1}(u_{k+1}-u_{k-2})\left[I+l_k(u_{k+1}-u_{k-1})\left[I+l_{k+1}(u_{k+1}-u_k)\right]^{-1}\right]^{-1}=$$

$$=\mathop{D}_{p=2}^{k}\frac{l_{k-p}\left(u_{k+1}-u_{k-p-1}\right)}{-I}=\mathop{D}_{p=1}^{k-1}\frac{l_{k-p-1}\left(u_{k+1}-u_{k-p-2}\right)}{-I},$$

. . . . .

$$I+l_2(u_{k+1}-u_1)[\cdots\left[I+l_{k-1}(u_{k+1}-u_{k-2})\times\right.$$

$$\left.\times\left[I+l_k(u_{k+1}-u_{k-1})\left[I+l_{k+1}(u_{k+1}-u_k)\right]^{-1}\right]^{-1}\right]^{-1}\cdots\Big]^{-1}=$$

$$=\left[I+l_0(u_{k+1})\right]^{-1}l_1(u_{k+1}-u_0)=\left[F(u_{k+1})-F(u_0)\right]^{-1}l_1(u_{k+1}-u_0),$$

that, according to the definition of abstract ICF (5), proves the validity of relation

$$T_n(u_{k+1}) = T_{k+1}(u_{k+1}) = F(u_{k+1}), \quad k = -1, 0, ..., n-1. \tag{9}$$

In a similar way the validity of continual interpolation conditions can be proved

$$T_n(u_{n-1,n}(\xi)) = F(u_{n-1,n}(\xi)), \quad \forall \xi \in [0,1]. \tag{10}$$

Thus, we have proved that the case is

**Theorem.** *Let us consider that* $F(u)$ *operator is* $n$ *times differentiated by Gateaux, and the one, for which fraction (5) makes sense. Then this fraction is abstract Thiel type ICF and satisfies the interpolation conditions (13), and continual interpolation condition (10).*

***Notice.*** Since abstract Thiel type ICF (5) satisfies the only one continual interpolation condition (10), its construction can be simplified by replacing the operator $g_\tau$ by operator $(1-\tau)I$. In such a case conditions (9) remain valid.

## 3. **Consequence 1**. Vector-valued and matrix-valued Thiele type ICF.

We are going to consider vector-valued and matrix-valued Thiel type interpolation. A significant number of works deal with vector-valued issue (see [1], [4], [15], [16] and referenced literature). Most of these works use the so-called Samelson inverse of vectors, which is a generally known procedure of matrix pseudoinverse of a full rank (see [14], 6.46). The same idea is actually used in the case of matrix-valued Thiel type interpolation, by previously turning matrices in the vectors (see [5], [17], [18] and referenced literature).

Let $Y$ algebra of $m \times m$ − matrices, $I = E$ a common identity matrix. Let us interpret the formulas the following way (5)-(7)

$$l_k(u-u_{k-1}) = -\frac{(u-u_{k-1})}{(u_k-u_{k-1})}\left[F_k(u_{k-1}) - F_k(u_k)\right],$$

$$k = 1, 2, ..., n, \quad F_1(u) = F(u),$$

$$F_i(u) = l_{i-1}\left(u-u_{i-2}\right)\left[-E + l_{i-2}\left(u-u_{i-3}\right)\times\right.$$

$$\left.\times\left[-E + l_{i-3}\left(u-u_{i-4}\right)\left[\cdots\left[-E - l_0(u)\right]^{-1}\cdots\right]^{-1}\right]^{-1}\right]^{-1} =$$

$$= \mathop{D}_{p=1}^{i} \frac{l_{i-p}\left(u-u_{i-p-1}\right)}{-E}, \ i = 1, 2, ..., u_{-1} = 0, \quad l_0(u) = F(u_0) - F(u) - E.$$

Consider a specific example for illustration.

**Example 1.** Let us consider that the following interpolation matrix conditions are defined

$$F(u_0)=\begin{bmatrix}2 & 0\\ 0 & -i\end{bmatrix},\quad F(u_1)=\begin{bmatrix}1 & 0\\ 1 & i\end{bmatrix},\quad F(u_2)=\begin{bmatrix}0 & i\\ 1 & 0\end{bmatrix},$$

$$u_0=-1,\ \ u_1=0,\ \ u_2=1.$$

Then

$$T_2(u)=\frac{1}{z^2-6z-3}\begin{bmatrix}-3-4z+7z^2 & -4i(z+1)z\\ -(z+1)(z+3) & i(-3+2z+z^2)\end{bmatrix},$$

which coincides with the example 2.8 з [18].

4. **Consequence 2. Matrix-valued Thiele type ICF for functionals from multiple variables.**

Let us consider the case of matrix-valued functionals from multiple variables. Here $X=R^k$, $Y$ - space of $m\times m$-matrices, which elements are functionals from $k$ variables, whose smoothness and range of definition are such that all the next formulas make sense. Hence, let

$$F(u)=\left[ f_{i,j}\underbrace{\left(\mathbf{x}(\cdot),\mathbf{y}(\cdot),\ldots,\mathbf{w}(\cdot)\right)}_{k} \right]_{i,j=1,\ldots,m}$$

and specified values of this matrix functional on continual knots

$$u=u_{s-1,s}(\tau)=\|x_{s-1}(t)+\tau(x_s(t)-x_{s-1}(t)),y_{s-1}(t)+\tau(y_s(t)-y_{s-1}(t)),\ldots,$$

$$\ldots,w_{s-1}(t)+\tau(w_s(t)-w_{s-1}(t))\|,\quad \tau\in[0,1],$$

$$F(u_{s-1,s}(\tau)),\quad s=1,2,\ldots,n.$$

Let us make formulas (6), (7) more specific for the consequence under consideration

$$l_r(u-u_{r-1})=\int_0^1 F_r{}'(u_{r-1,r}(\tau_r))(u-u_{r-1})d\tau_r=$$

$$=\int_0^1 F_r{}'(u_{r-1}+\tau_r(u_r-u_{r-1}))(u-u_{r-1})d\tau_r,$$

$$r=1,2,\ldots,n, \quad F_1(u)=F(u),$$

$$F_r'(u_{r-1,r}(\tau_r))(u-u_{r-1})=\left[\frac{\partial}{\partial x} f_{i,j}(u_{r-1,r}(\tau_r))\right]_{i,j=1,2,\ldots,m}\left(x(\cdot)-x_{r-1}(\cdot)\right)+$$

$$+\left[\frac{\partial}{\partial y} f_{i,j}(u_{r-1,r}(\tau_r))\right]_{i,j=1,2,\ldots,m}\left(y(\cdot)-y_{r-1}(\cdot)\right)+\ldots$$

$$\ldots+\left[\frac{\partial}{\partial w} f_{i,j}(u_{r-1,r}(\tau_r))\right]_{i,j=1,2,\ldots,m}\left(w(\cdot)-w_{r-1}(\cdot)\right),$$

$$F_i(u)=\sum_{p=1}^{i}\frac{l_{i-p}\left(u-u_{i-p-1}\right)}{-I}, \quad i=1,2,\ldots,$$

$$u_{-1}=0, \quad l_0(u)=F(u_0)-F(u)-I.$$

**Example 2.** Let $X=R^2$ and

$$F(u)=\left[f_{i,j}(x,y)\right]_{i,j=1,2}=\begin{bmatrix}\sin(x+y), \cos(x+y)\\ x^2 \quad , 1/(1+y)\end{bmatrix}.$$

Then, given such interpolation knots

$$u_0=(0,0), \quad u_1=(\pi/2,\pi/2), \quad u_2=(\pi,0)$$

we get

$$F(u_0+t(u_1-u_0))=\left[f_{i,j}(t\pi/2,t\pi/2)\right]_{i,j=0,1},$$

$$l_1(u-u_0)=\begin{bmatrix}0 & ,-\frac{2}{\pi}(x+y)\\ \frac{\pi}{2}x, & -\frac{2y}{2+\pi}\end{bmatrix},$$

$$l_2(u-u_1)=\int_0^1 F_2'(u_1+\tau_2(u_2-u_1))(u-u_1)d\tau_2,$$

$$F_2(u)=l_1(u)\left[F(u)-F(u_0)\right]^{-1}=$$

$$=\Delta^{-1}(x,y)\begin{bmatrix}0 & -\frac{2}{\pi}(x+y)\\ \frac{\pi}{2}x & -\frac{2y}{2+\pi}\end{bmatrix}\begin{bmatrix}-\frac{y}{1+y} & 1-\cos(x+y)\\ -x^2 & \sin(x+y)\end{bmatrix}=$$

$$=\Delta^{-1}(x,y)\begin{bmatrix}\frac{2}{\pi}(x+y)x^2 & -\frac{2}{\pi}(x+y)\sin(x+y)\\ -\frac{\pi xy}{2(1+y)}+\frac{2yx^2}{2+\pi} & \frac{\pi x}{2}\left(1-\cos(x+y)\right)-\frac{2y}{2+\pi}\sin(x+y)\end{bmatrix},$$

$$\Delta(x,y)=-\frac{y}{1+y}\sin(x+y)-x^2[\cos(x+y)-1],$$

$$l_2(u-u_1)=\begin{bmatrix}a_{1,1} & \frac{2(x+y-\pi)}{\pi^2}\\ a_{2,1} & a_{2,2}\end{bmatrix},$$

$$a_{1,1}=\frac{1}{\pi(\pi+1)^2}\left(\ln(2+\pi)+\pi+\pi^2\right)(x+y-\pi),$$

$$a_{2,1}=\frac{1}{8\pi(2+\pi)(1+\pi)^4}\Big(8(1+\pi)^4\ln 2-2\left(2\pi^4+10\pi^3+17\pi^2+10\pi+4\right)\ln(2+\pi)+$$
$$+\pi(1+\pi)\left(4\pi^3+11\pi^2+11\pi-2\right)\Big)(x+y-\pi),$$

$$a_{2,2}=-\frac{1}{4\pi(2+\pi)(1+\pi)^3}\Big[(1+\pi)\Big(\left(4\pi^3+11\pi^2+11\pi+2\right)(x-\pi/2)-$$
$$-\left(5\pi^2+9\pi+6\right)(y-\pi/2)\Big)+$$
$$+\left(8(1+\pi)^3\ln 2-2\pi(2+\pi)\ln(2+\pi)\right)(x+y-\pi/2)\Big]$$

and, as a result

$$T_2(x,y)=$$
$$=\frac{1}{\Delta_2(x,y)}\begin{bmatrix}t_{1,1}(\mathrm{x,y}) & t_{1,2}(\mathrm{x,y})\\ t_{2,1}(\mathrm{x,y}) & t_{2,2}(x,y)\end{bmatrix},$$

$$\Delta_2(x,y)=0.10094x^2-0.45299x+0.11535xy-0.07314-0.49946y+0.01441y^2,$$

$$t_{1,1}(\mathrm{x},\mathrm{y}) = 1.5708(0.20264x - 0.63662 + 0.20264y)x,$$

$$\mathrm{t}_{1,2}(\mathrm{x},\mathrm{y}) = -0.10097x^2 + 0.49963x - 0.16465xy - 0.07314 + 0.7008y - 0.06368y^2,$$

$$t_{2,1}(x,y) = -1.5708(0.27184x + 0.146 + 0.27184y)x,$$

$$t_{2,2}(x,y) = 0.05468x^2 - 0.30766x + 0.12857xy - 0.07314 - 0.29734y + 0.07389y^2.$$